\documentclass{amsart}
\usepackage{setspace}

\usepackage{amsmath,amsthm,amsfonts,amscd,amssymb}

\newtheorem{thm}{Theorem}[section]
\newtheorem{cor}[thm]{Corollary}
\newtheorem{lem}[thm]{Lemma}

\theoremstyle{definition}

\theoremstyle{remark}
\newtheorem{rem}{Remark}[section]

\begin{document}

\title{On well-rounded ideal lattices}
\author{Lenny Fukshansky and Kathleen Petersen}

\address{Department of Mathematics, 850 Columbia Avenue, Claremont McKenna College, Claremont, CA 91711}
\email{lenny@cmc.edu}
\address{Department of Mathematics, Florida State University, 208 Love Building, 1017 Academic Way, Tallahassee, FL 32306}
\email{petersen@math.fsu.edu}

\begin{abstract}
We investigate a connection between two important classes of Euclidean lattices: well-rounded and ideal lattices. A lattice of full rank in a Euclidean space is called well-rounded if its set of minimal vectors spans the whole space. We consider lattices coming from full rings of integers in number fields, proving that only cyclotomic fields give rise to well-rounded lattices. We further study the well-rounded lattices coming from ideals in quadratic rings of integers, showing that there exist infinitely many real and imaginary quadratic number fields containing ideals which give rise to well-rounded lattices in the plane.
\end{abstract}

\maketitle

\def\A{{\mathcal A}}
\def\AA{{\mathfrak A}}
\def\B{{\mathcal B}}
\def\C{{\mathcal C}}
\def\D{{\mathcal D}}
\def\EE{{\mathfrak E}}
\def\F{{\mathcal F}}
\def\x{{\mathcal H}}
\def\I{{\mathcal I}}
\def\II{{\mathfrak I}}
\def\J{{\mathcal J}}
\def\K{{\mathcal K}}
\def\kk{{\mathfrak K}}
\def\L{{\mathcal L}}
\def\LL{{\mathfrak L}}
\def\M{{\mathcal M}}
\def\mm{{\mathfrak m}}
\def\MM{{\mathfrak M}}
\def\N{{\mathcal N}}
\def\O{{\mathcal O}}
\def\OO{{\mathfrak O}}
\def\PP{{\mathfrak P}}
\def\R{{\mathcal R}}
\def\PNR{{\mathcal P_N(\real)}}
\def\PMNR{{\mathcal P^M_N(\real)}}
\def\PdNR{{\mathcal P^d_N(\real)}}
\def\s{{\mathcal S}}
\def\V{{\mathcal V}}
\def\X{{\mathcal X}}
\def\Y{{\mathcal Y}}
\def\Z{{\mathcal Z}}
\def\H{{\mathcal H}}
\def\cee{{\mathbb C}}
\def\pee{{\mathbb P}}
\def\que{{\mathbb Q}}
\def\QQ{{\mathbb Q}}
\def\real{{\mathbb R}}
\def\RR{{\mathbb R}}
\def\zed{{\mathbb Z}}
\def\ZZ{{\mathbb Z}}
\def\aaa{{\mathbb A}}
\def\ff{{\mathbb F}}
\def\kk{{\mathfrak K}}
\def\qbar{{\overline{\mathbb Q}}}
\def\kbar{{\overline{K}}}
\def\ybar{{\overline{Y}}}
\def\kkbar{{\overline{\mathfrak K}}}
\def\ubar{{\overline{U}}}
\def\eps{{\varepsilon}}
\def\ahat{{\hat \alpha}}
\def\bhat{{\hat \beta}}
\def\gt{{\tilde \gamma}}
\def\h{{\tfrac12}}
\def\be{{\boldsymbol e}}
\def\bei{{\boldsymbol e_i}}
\def\bc{{\boldsymbol c}}
\def\bm{{\boldsymbol m}}
\def\bk{{\boldsymbol k}}
\def\bi{{\boldsymbol i}}
\def\bl{{\boldsymbol l}}
\def\bq{{\boldsymbol q}}
\def\bu{{\boldsymbol u}}
\def\bt{{\boldsymbol t}}
\def\bs{{\boldsymbol s}}
\def\bv{{\boldsymbol v}}
\def\bw{{\boldsymbol w}}
\def\bx{{\boldsymbol x}}
\def\bX{{\boldsymbol X}}
\def\bz{{\boldsymbol z}}
\def\bZ{{\boldsymbol Z}}
\def\bwy{{\boldsymbol y}}
\def\bY{{\boldsymbol Y}}
\def\bL{{\boldsymbol L}}
\def\ba{{\boldsymbol a}}
\def\bb{{\boldsymbol b}}
\def\bet{{\boldsymbol\eta}}
\def\bxi{{\boldsymbol\xi}}
\def\balpha{{\boldsymbol\alpha}}
\def\bo{{\boldsymbol 0}}
\def\bone{{\boldsymbol 1}}
\def\bol{{\boldsymbol 1}_L}
\def\ep{\varepsilon}
\def\p{\boldsymbol\varphi}
\def\q{\boldsymbol\psi}
\def\rank{\operatorname{rank}}
\def\aut{\operatorname{Aut}}
\def\lcm{\operatorname{lcm}}
\def\sgn{\operatorname{sgn}}
\def\spn{\operatorname{span}}
\def\md{\operatorname{mod}}
\def\Norm{\operatorname{Norm}}
\def\dim{\operatorname{dim}}
\def\det{\operatorname{det}}
\def\Vol{\operatorname{Vol}}
\def\rk{\operatorname{rk}}
\def\ord{\operatorname{ord}}
\def\ker{\operatorname{ker}}
\def\div{\operatorname{div}}
\def\Gal{\operatorname{Gal}}
\def\GL{\operatorname{GL}}
\def\SNR{\operatorname{SNR}}
\def\WR{\operatorname{WR}}
\def\Disc{\operatorname{Disc}}
\def\scg{\operatorname{\left< \Gamma \right>}}
\def\swrh{\operatorname{Sim_{WR}(\Lambda_h)}}
\def\ch{\operatorname{C_h}}
\def\cht{\operatorname{C_h(\theta)}}
\def\scgt{\operatorname{\left< \Gamma_{\theta} \right>}}
\def\scgmn{\operatorname{\left< \Gamma_{m,n} \right>}}
\def\gat{\operatorname{\Omega_{\theta}}}
\def\NN{\operatorname{N}}

\section{Introduction}
\label{intro}

In this note we investigate a connection between two fundamental classes of Euclidean lattices, {\it well-rounded} and {\it ideal} lattices, which come up in a variety of mathematical contexts as well as in applications in discrete optimization and coding theory.
\smallskip

Let $\Lambda$ be a lattice of full rank in the $d$-dimensional Euclidean space $\real^d$ for $d \geq 2$. The {\it minimum} of $\Lambda$ is defined as
$$|\Lambda| := \min \{ \| \bx \|^2 : \bx \in \Lambda \setminus \{ \bo \} \},$$
where $\|\ \|$ stands for the usual Euclidean norm on $\real^d$, and the set of {\it minimal vectors} of $\Lambda$ is defined to be
$$S(\Lambda) := \{ \boldsymbol x \in \Lambda :  \| \boldsymbol x \|^2 = |\Lambda| \}.$$
The lattice $\Lambda$ is called {\it well-rounded} (abbreviated WR) if the set $S(\Lambda)$ spans $\real^d$. WR lattices are important in discrete optimization, in particular in the investigation of sphere packing, sphere covering, and kissing number problems (see \cite{martinet}), as well as in coding theory (see \cite{esm}). Properties of WR lattices have also been investigated in \cite{mcmullen} in connection with Minkowski's conjecture.

Another class of lattices that comes up frequently in connection with optimization problems and in coding theory (see \cite{tsfasman}, \cite{oggier}) are the ideal lattices. Let $K$ be a number field of degree $d$ over $\que$, and let us write $\O_K$ for its ring of integers. Let 
$$\sigma_1,...,\sigma_{r_1},\tau_1,...,\tau_{r_2},...,\tau_{2r_2}$$
be the embeddings of $K$ into $\cee$ with $\sigma_1,...,\sigma_{r_1}$ being the real embeddings and $\tau_n,\tau_{r_2+n} = \bar{\tau}_n$ for each $1 \leq n \leq r_2$ being the pairs of complex conjugate embeddings. For each $\alpha \in K$ and each complex embedding $\tau_n$, write $\tau_{n1}(\alpha) = \Re(\tau_n(\alpha))$ and $\tau_{n2}(\alpha) = \Im(\tau_n(\alpha))$, where $\Re$ and $\Im$ stand respectively for real and imaginary parts of a complex number. Then $d=r_1+2r_2$, and we define an embedding
$$\sigma = (\sigma_1,\dots,\sigma_{r_1}, \tau_{11},\tau_{12},\dots,\tau_{r_21},\tau_{r_22}): K \to \real^d.$$
Then $\Lambda_K := \sigma(\O_K)$ is a lattice of full rank in $\real^d$. Following the notation of \cite{bayer_nebe} (bottom of p. 438), we call such lattices {\it principal ideal lattices}. More generally, for any nonzero fractional ideal $I$ of $\O_K$, $\Lambda_K(I) := \sigma(I)$ is a full rank lattice in $\real^d$, and if $I$ is an ideal in $\O_K$ then $\Lambda_K(I)$ is a sublattice of $\Lambda_K$ of finite index; throughout this paper, when we refer to ideals or fractional ideals, we always mean only the nonzero ones. A lattice $\Lambda$ in $\real^d$ is called an {\it ideal lattice} if it can be realized as $\Lambda_K(I)$ for some fractional ideal $I$ of the ring of integers of some number field $K$ with $[K:\que]=d$. For more information on ideal lattices see \cite{bayer1}, \cite{bayer2}, \cite{bayer_nebe}. It should be remarked that the definition of ideal lattices (and principal ideal lattices in particular) in these papers is more general, our definition being a more concrete special case of that.

The importance and applicability of these two special classes of lattices motivates the following natural question: when are ideal lattices well-rounded?  In this note we investigate this question for principal ideal lattices $\Lambda_K$ and some of their ideal sublattices. This question is partially motivated by the first author's previous investigations \cite{wr1}, \cite{wr2}, \cite{wr3}, where the WR sublattices of $\zed^2 = \Lambda_{\que(i)}$ and the hexagonal lattice $\Lambda_h := \Lambda_{\que(\sqrt{-3})}$ were studied. Both of these lattices are WR themselves; in fact, these are the only two principal ideal WR lattices in $\real^2$, as we demonstrate in Section~\ref{quadratic} by a direct verification argument (see Lemma~\ref{O_K_WR}). Moreover, all ideal sublattices of $\zed^2$ and $\Lambda_h$ are also WR: this is a direct consequence of the well-known fact the ideal sublattices of $\zed^2$ and $\Lambda_h$ are similar to $\zed^2$ and $\Lambda_h$, respectively. Let us recall here that two lattices $\Lambda$ and $\Omega$ are said to be {\it similar} if there exists an $N \times N$ real orthogonal matrix $A$ and a nonzero constant $\alpha$ such that $\Omega = \alpha A \Lambda$. Similarity is easily seen to be an equivalence relation, which preserves the WR property; we will denote it by writing $\Lambda \sim \Omega$. In Section~\ref{quadratic} we further investigate the ideal lattices coming from quadratic number fields, proving in particular the following result.

\begin{thm} \label{q_ideals} There exist infinitely many real and imaginary quadratic number fields $K$ whose rings of integers contain an ideal $I$ such that the planar lattice $\Lambda_K(I)$ is WR.
\end{thm}

\noindent
We give examples of ideals as in Theorem \ref{q_ideals} in Tables~\ref{table1} and~\ref{table2} in Section~\ref{quadratic}. 
\smallskip

\begin{rem} \label{qi-qr_ideals} We should remark that there also exist quadratic number fields $K$ so that $\Lambda_K(I)$ is not WR for any ideal $I \subseteq \O_K$, for instance all class number one imaginary quadratic fields different from $\que(i)$ and $\que(\sqrt{-3})$, as we demonstrate in Corollary~\ref{cycl_pid}.
\end{rem}

The distinguishing feature of $\que(i)$ and $\que(\sqrt{-3})$ ($=\que(e^{\frac{2\pi i}{3}})$) among imaginary quadratic number fields is that these are the only ones that are cyclotomic fields. In Section~\ref{cyclotom} we show that for number fields of any degree principal ideal lattices are WR only in the cyclotomic case. Recall that the {\it norm} of a nonzero ideal $I \subseteq \O_K$ in the number field $K$ is defined as $\NN(I) := |\O_K/I|$. We prove the following result.

\begin{thm} \label{main} Let $K$ be a number field of degree $d \geq 2$ and $I \subseteq \O_K$ a nonzero ideal. Then $|\Lambda_K(I)| \geq (r_1+r_2) \NN(I)^{\frac{1}{r_1+r_2}}$. Moreover, $|\Lambda_K|=r_1+r_2$,
$$S(\Lambda_K) = \{ \sigma(x) : x \in \O_K \text{ is a root of unity} \},$$
and $\Lambda_K$ is WR if and only if $K$ is a cyclotomic field, i.e., $K=\que(\zeta_k)$ for some primitive $k$-th root of unity $\zeta_k$, $k \geq 2$. If this is the case, then 
$$|\Lambda_K| = r_2 = \frac{d}{2} =  \frac{\varphi(k)}{2}.$$
\end{thm}
\smallskip

\begin{rem} \label{bayer} The value of $|\Lambda_K|$ and a bound on $|\Lambda_K(I)|$ under a slightly different embedding into $\real^d$ also follow from Lemma 4.3 of \cite{bayer_eucl}, which is proved by a rather different argument from ours, however the results of \cite{bayer_eucl} do not imply the result of Theorem~\ref{main} on WR ideal lattices.
\end{rem}

\noindent
As a corollary of Theorem~\ref{main}, we also deduce that, as in the two dimensional case, all ideal  lattices coming from cyclotomic fields are WR.

\begin{cor} \label{cycl_ideal} Let $K=\que(\zeta_k)$ for some primitive $k$-th root of unity $\zeta_k$, $k \geq 2$, and let $I$ be a fractional ideal of $\O_K$. Then the lattice $\Lambda_K(I)$ is WR.
\end{cor}

\noindent
The proof of Corollary~\ref{cycl_ideal} is also presented in Section~\ref{cyclotom}. 
\bigskip

\section{Quadratic ideal lattices}
\label{quadratic}

In this section we study quadratic WR ideal lattices. Let us start by recording a general basic property of WR lattices in $\real^2$ which will be useful to us.

\begin{lem} \label{pairs} A full-rank lattice $\Lambda \subset \real^2$ contains 2, 4, or 6 minimal vectors, and it is WR if and only if $|S(\Lambda)| = 4,6$. Moreover, $|S(\Lambda)|=6$ if and only if $\Lambda$ is similar to $\Lambda_h$, the hexagonal lattice. On the other hand, there are infinitely many distinct similarity classes of WR lattices in $\real^2$ with four minimal vectors.
\end{lem}

\proof
Notice that minimal vectors in a lattice always come in $\pm$ pairs, hence $S(\Lambda)$ contains at least two vectors, and it contains two linearly independent vectors if and only if its cardinality is greater than two. On the other hand, the angle $\theta$ between any pair of minimal vectors $\bx,\bwy$ has to be at least $\pi/3$, since otherwise
$$\|\bx-\bwy\|^2 = \|\bx\|^2 + \|\bwy\|^2 - 2\|\bx\| \|\bwy\| \cos \theta < \|\bx\|^2 = \|\bwy\|^2.$$
Since all the minimal vectors must lie on the circle of the same radius,
$$|S(\Lambda)| \leq \frac{2\pi}{\pi/3} = 6,$$
and so $\Lambda$ is WR if and only if $|S(\Lambda)| = 4,6$.

Now suppose that $|S(\Lambda)| = 6$, then 
$$S(\Lambda) = \{ \pm \bx_1, \pm \bx_2, \pm \bx_3 \},$$
and we can choose a pair $\pm \bx_i,\pm \bx_j$, $1 \leq i < j \leq 3$, such that the angle between these two vectors is $\pi/3$. Then $\Lambda$ is spanned over $\zed$ by these two vectors, and hence can be obtained from $\Lambda_h$ by rotation and dilation, i.e. is similar to $\Lambda_h$. On the other hand, if $\Lambda$ is similar to $\Lambda_h$, then
$$|S(\Lambda)| = |S(\Lambda_h)| = 6.$$

On the other hand, there are infinitely many distinct similarity classes of WR lattices in $\real^2$ (see \cite{wr2}), and so, by our argument above, lattices in all but one of them cannot contain six minimal vectors; hence they must have four. This completes the proof.
\endproof

\begin{rem} \label{pairs_known} The statement of Lemma~\ref{pairs} is generally well-known with parts of it following from the work of Gauss (see Section~3 of \cite{wr1} for some details). We present the proof here for completeness purposes.
\end{rem}

\smallskip

We first consider principal ideal lattices in $\real^2$. Let $K = \que(\sqrt{D})$ for some squarefree $D \in \zed$, $D \neq 1$, then
\begin{equation}
\label{ring_int}
\O_K = \left\{ \begin{array}{ll}
\zed[\sqrt{D}] & \mbox{if $D \not\equiv 1 (\md 4)$} \\
\zed[\frac{1+\sqrt{D}}{2}] & \mbox{if $D \equiv 1 (\md 4)$,}
\end{array}
\right.
\end{equation}
and the two embeddings of $K$ are given by
$$\sqrt{D} \mapsto \sqrt{D},\ \sqrt{D} \mapsto -\sqrt{D}.$$
Let us use the notation $\Lambda_D$ for $\Lambda_K$ and $\Lambda_D(I)$ for any ideal $I \subset \O_{\que(\sqrt{D})}$, so for instance $\zed^2 = \Lambda_{-1}$ and $\Lambda_h = \Lambda_{-3}$. Our first lemma stipulates that these are the only cases when $\Lambda_D$ is well-rounded. While it is a special case of Lemma~\ref{only_cyclotomic} below, we prove it here by a direct elementary argument.

\begin{lem} \label{O_K_WR} The lattice $\Lambda_D$ is WR if and only if $D=-1,-3$.
\end{lem}

\proof First assume that $D \not\equiv 1 (\md 4)$ is positive, then
$$\Lambda_D = \left( \begin{matrix} 1&\sqrt{D} \\ 1&-\sqrt{D} \end{matrix} \right) \zed^2.$$
Now for any nonzero
$$\bx = \left( \begin{matrix} 1&\sqrt{D} \\ 1&-\sqrt{D} \end{matrix} \right) \left( \begin{matrix} m \\ n \end{matrix} \right) = \left( \begin{matrix} m+n\sqrt{D} \\ m-n\sqrt{D} \end{matrix} \right) \in \Lambda_D$$
we have 
$$\|\bx\|^2 = (m+n\sqrt{D})^2 + (m-n\sqrt{D})^2 = 2(m^2+Dn^2) \geq 2,$$
with equality in this inequality if and only if $m = \pm 1, n = 0$, which means that $|\Lambda_D| = 2$ and
$$S(\Lambda_D) = \left\{ \pm \left( \begin{matrix} 1 \\ 1 \end{matrix} \right) \right\},$$
and so $\Lambda_D$ cannot be WR.
\smallskip

Next assume that $D \not\equiv 1 (\md 4)$ is negative, then
$$\Lambda_D = \left( \begin{matrix} 1&0 \\ 0&\sqrt{|D|} \end{matrix} \right) \zed^2.$$
Hence $|\Lambda_D| = 1$, and for any nonzero
$$\bx = \left( \begin{matrix} m \\ n\sqrt{|D|} \end{matrix} \right) \in \Lambda_D$$
we have $\|\bx\|^2 = m^2 + |D|n^2 \geq 1$, with equality in this inequality if and only if $m = \pm 1, n = 0$, unless $D=-1$, in which case there are additional solutions $m = 0, n = \pm 1$. Hence $\Lambda_D$ is not WR unless $D=-1$, and in this later case
$$S(\Lambda_{-1}) = \left\{ \pm \left( \begin{matrix} 1 \\ 0 \end{matrix} \right),   \pm \left( \begin{matrix} 0 \\ 1 \end{matrix} \right) \right\}.$$
Therefore for negative $D \not\equiv 1 (\md 4)$, $\Lambda_D$ is WR if and only if $D=-1$.
\smallskip

Next assume that $D \equiv 1 (\md 4)$ is positive, then $D \geq 5$ and
$$\Lambda_D = \left( \begin{matrix} 1&\frac{1+\sqrt{D}}{2} \\ 1&\frac{1-\sqrt{D}}{2} \end{matrix} \right) \zed^2,$$
and so for any nonzero
$$\bx = \left( \begin{matrix} \frac{2m+n}{2} + \frac{n\sqrt{D}}{2} \\ \frac{2m+n}{2} - \frac{n\sqrt{D}}{2} \end{matrix} \right) \in \Lambda_D$$
we have 
$$\|\bx\|^2 = \frac{1}{2} \left( 4m^2 + (D+1)n^2 + 4mn \right) \geq 2 m^2 + 3n^2 + 2mn \geq 2,$$
with equality in this inequality if and only if $m = \pm 1, n = 0$, which means that $|\Lambda_D| = 2$ and
$$S(\Lambda_D) = \left\{ \pm \left( \begin{matrix} 1 \\ 1 \end{matrix} \right) \right\},$$
and so $\Lambda_D$ cannot be WR.
\smallskip

Finally suppose that $D \equiv 1 (\md 4)$ is negative, then $D \leq -3$ and
$$\Lambda_D = \left( \begin{matrix} 1&\frac{1}{2} \\ 0&\frac{\sqrt{|D|}}{2} \end{matrix} \right) \zed^2.$$
Hence $|\Lambda_D| = 1$, and for any nonzero
$$\bx = \left( \begin{matrix} \frac{2m+n}{2} \\ \frac{n\sqrt{|D|}}{2} \end{matrix} \right) \in \Lambda_D$$
we have
$$\|\bx\|^2 = m^2 + mn + \frac{(|D|+1) n^2}{4} \geq  1$$
with equality if and only if $m = \pm 1, n = 0$, unless $D=-3$, in which case there are additional solutions $m=0, n=\pm 1$, and $m=1, n= -1$, as well as $m=-1, n= 1$. Hence $\Lambda_D$ is not WR unless $D=-3$, and in this later case
$$S(\Lambda_{-3}) = \left\{ \pm \left( \begin{matrix} 1 \\ 0 \end{matrix} \right),   \pm \left( \begin{matrix} \frac{1}{2} \\ \frac{\sqrt{3}}{2} \end{matrix} \right), \pm \left( \begin{matrix} \frac{1}{2} \\ -\frac{\sqrt{3}}{2} \end{matrix} \right) \right\}.$$
Therefore for negative $D \equiv 1 (\md 4)$, $\Lambda_D$ is WR if and only if $D=-3$. This completes the proof.
\endproof
\bigskip

Next we discuss more general WR ideal lattices coming from quadratic number fields. Notice that if $K$ is a quadratic number field, then either it is an imaginary quadratic field (i.e. $K=\que(\sqrt{D})$ with $D \leq -1$ a squarefree integer, so that $r_1=0,r_2=1$) or a real quadratic field (i.e. $K=\que(\sqrt{D})$ with $D > 1$ a squarefree integer, so that $r_1=2,r_2=0$). We first consider imaginary quadratic fields, and start by establishing a basic property of principal ideals.

\begin{lem}\label{lemma:principalsimilarity} Let $K$ be an imaginary quadratic number field. If $I \subseteq \O_K$ is a principal ideal and $J=\alpha I$, $0 \neq \alpha \in K$, is a fractional ideal, then $\Lambda_K(J)$ is similar to $\Lambda_K$.
\end{lem}

\proof
Since $I$ is a principal ideal, $I=\gamma \O_K$ for some $\gamma \in \O_K$, and so $J=\alpha' \O_K$, where $\alpha'= \alpha \gamma \in \cee$.  Writing $\alpha' = re^{i\theta}$ for $r, \theta \in \real$, the action of left multiplication by $\alpha'$ on an element $\beta = se^{i\phi}$ is $\alpha'  \beta = rse^{i(\theta +\phi)}$ which is a dilation and a rotation.  Since $\Lambda_K = \sigma(\O_K)$, this is the action of $\alpha'$ on the lattice, meaning that $\Lambda_K(J)$ is obtained from $\Lambda_K$ by rotation and dilation. Hence the two lattices are similar.
\endproof

\begin{cor} \label{cycl_pid} Let $K=\que(i)$ or $\que(\sqrt{-3})$ and $I \subseteq K$ a fractional ideal, then $\Lambda_K(I)$ is WR. On the other hand, if $K$ is an imaginary quadratic field $ \neq \que(i), \que(\sqrt{-3})$ and $I$ is a principal fractional ideal in $K$, then $\Lambda_K(I)$ is not WR.
\end{cor}

\proof
Both of the fields $\que(i)$ or $\que(\sqrt{-3})$ are principal ideal domains, and so the first statement follows by combining Lemma~\ref{O_K_WR} with Lemma~\ref{lemma:principalsimilarity}. On the other hand, $\Lambda_K$ is not WR whenever $K \neq \que(i), \que(\sqrt{-3})$ by Lemma~\ref{O_K_WR}, and so the second statement follows by Lemma~\ref{lemma:principalsimilarity}.
\endproof

We will next construct an infinite family of imaginary quadratic fields with ideals giving rise to WR ideal lattices. Our construction is based on a certain convenient choice of an integral basis for an ideal in any quadratic number field. Let $D$ be a squarefree integer, and define
\begin{equation}
\label{delta_def}
\delta = \left\{ \begin{array}{ll}
-\sqrt{D} & \mbox{if $D \not\equiv 1 (\md 4)$} \\
\frac{1-\sqrt{D}}{2} & \mbox{if $D \equiv 1 (\md 4)$.}
\end{array}
\right.
\end{equation}
Let $K=\que(\sqrt{D})$ and let $I$ be an ideal in $\O_K$, where $\O_K$ is as in \eqref{ring_int}. It is a well-known fact (see for instance Theorem 6.9 on p. 94 of \cite{buell}) that there exist rational integers $a,b,g$ with
\begin{equation}
\label{abg}
0 \leq b < a,\ 0 < g \leq a,\ g \mid a,\ g \mid b,
\end{equation}
so that
\begin{equation}
\label{I_abg}
I = \{ ax + (b+g\delta)y : x,y \in \zed \}.
\end{equation}
In other words, $a, b+g\delta$ is an integral basis for $I$; moreover, an integral basis for $I$ with these properties is unique. Further, if a triple $a,b,g \in \zed$ satisfying \eqref{abg} in addition satisfies the condition
\begin{equation}
\label{norm_abg}
\NN(b+g\delta) = kga, \text{ for some integer } k,
\end{equation}
where $\NN$ stands for the norm, then the corresponding ideal $I = \left< a, b+g\delta \right>$ in $\O_K$ is of the form \eqref{I_abg} (see Theorem 6.15 on p. 96 of \cite{buell}). The unique integral basis with these properties is called the canonical basis for the ideal. Our strategy in the arguments to follow is based on using the canonical basis for an ideal $I$ in $\O_K$ to construct a basis for the lattice $\Lambda_K(I)$ whose corresponding norm form is Minkowski reduced. Given a basis matrix $A = (\bx\ \bwy)$ for a lattice in $\real^2$, the corresponding norm form is 
$$Q(m,n) = c_1 m^2 + c_2 mn + c_3 n^2 := (m\ n) A^t A \left( \begin{matrix} m \\ n \end{matrix} \right).$$ 
It is said to be reduced if $Q(m,n) \geq Q(0,1) \geq Q(1,0)$ for all $m,n \in \zed$ with $n \neq 0$, which is equivalent to saying that $|c_2| \leq c_1 \leq c_3$. If in addition $Q$ is symmetric, meaning that $Q(1,0) = Q(0,1)$ (i.e., $c_1=c_3$), the lattice must be WR.
\smallskip

We return to imaginary quadratic fields.

\begin{lem} \label{im_quad} There exist infinitely many squarefree integers $D > 1$ with $-D \equiv 1 (\md 4)$ for which the ring of integers $\O_K$ of the imaginary quadratic field $K=\que(\sqrt{-D})$ contains an ideal $I$ with the property that $\Lambda_K(I)$ is WR.
\end{lem}

\proof
Let $t$ be an odd positive integer, and define
\begin{eqnarray}
\label{D_ab}
& \ & g=1,\ b=\frac{t-1}{2},\ a=2b+2=t+1, \nonumber \\
& \ & D = (t+2)(3t+2) = 3t^2+8t+4.
\end{eqnarray}
It is then easy to see that $D \equiv 3 (\md 4)$, i.e.  $-D \equiv 1 (\md 4)$, and condition \eqref{abg} is satisfied. Moreover, there exist infinitely many odd positive $t$ for which $D$ given by \eqref{D_ab} is squarefree. Indeed, notice that the set $\{ t+2 : t \in \zed,\ 2 \nmid t \}$ contains the set $\pee$ of all odd prime numbers. Then select $t$ such that $p=t+2 \in \pee$, so
$$D=p(3p-4),$$
and clearly $p \nmid 3p-4$. In this case, to ensure that $D$ is squarefree we only need to select $t$ in such a way that $3p-4$ is squarefree. The fact that there exist infinitely many prime numbers $p$ such that $3p-4$ is squarefree follows from the theorem on p. 920 of \cite{clary}; for each such prime $p$, let $t=p-2$. For each such choice of $t$, let $K=\que(\sqrt{-D})$ and 
$$I = \left< 2b+2, b+\delta \right>  = \left< t+1, \frac{t}{2} - \frac{\sqrt{-D}}{2} \right> \subseteq \O_K$$
be an ideal. Then
\begin{eqnarray*}
\NN(b+\delta) & = & \left(\frac{t-\sqrt{-D}}{2} \right) \left(\frac{t+\sqrt{-D}}{2} \right) = \frac{1}{4} (t^2 + D) \\
& = & \frac{1}{4} (4t^2+8t+4) = (t+1)^2 = a^2,
\end{eqnarray*}
and so the condition \eqref{norm_abg} is satisfied. Therefore $t+1, \frac{t}{2} - \frac{\sqrt{-D}}{2}$ is a canonical basis for $I$. Then $\Lambda_K(I) = A\zed^2$, where
\begin{equation}
\label{A_1}
A = \left( \begin{matrix} t+1 & \frac{t}{2} \\ 0 & - \frac{\sqrt{D}}{2} \end{matrix} \right).
\end{equation}
Then for any $\bx = A \left( \begin{matrix} m \\ n \end{matrix} \right) \in \Lambda_K(I)$,
\begin{eqnarray}
\label{Q_mn}
\|\bx\|^2 & = & Q(m,n) := (m\ n) A^t A \left( \begin{matrix} m \\ n \end{matrix} \right) \nonumber \\
& = & (t+1)^2 m^2 + t(t+1) mn + \frac{1}{4}(t^2+D) n^2 = a^2 m^2 + a(a-1) mn + a^2 n^2,
\end{eqnarray}
and hence the positive definite integral binary quadratic form $Q(m,n)$ is reduced and symmetric. By definition of Minkowski reduction, $Q(m,n) \geq Q(0,1) \geq Q(1,0)$ for all $m,n \in \zed$ with $n \neq 0$, and by symmetry $Q(0,1)=Q(1,0)=a^2$. This implies that $\Lambda_K(I)$ is WR, and thus we have constructed an infinite family of imaginary quadratic number fields $\que(\sqrt{-D})$ with $D>1$, $-D \equiv 1 (\md 4)$, each of which contains at least one ideal $I \subseteq \O_K$ so that $\Lambda_K(I)$ is WR. This completes the proof of the lemma.
\endproof
\smallskip

We now turn to the case of real quadratic number fields, and establish a result analogous to Lemma~\ref{im_quad}.

\begin{lem} \label{real_quad} There exist infinitely many squarefree integers $D > 1$ with $D \equiv 1 (\md 4)$ for which the ring of integers $\O_K$ of the real quadratic field $K=\que(\sqrt{D})$ contains an ideal $I$ with the property that $\Lambda_K(I)$ is WR.
\end{lem}

\proof
Let $t$ be an odd positive integer, and define
\begin{eqnarray}
\label{D_ab_r}
& \ & g=1,\ b=\frac{t+1}{2},\ a=2b+1=t+2, \nonumber \\
& \ & D = (t+2)(t-2) = t^2-4.
\end{eqnarray}
It is then easy to see that $D \equiv 1 (\md 4)$ and condition \eqref{abg} is satisfied. Moreover, we can show that there exist infinitely many odd positive $t$ for which $D$ given by \eqref{D_ab_r} is squarefree, using the same type of argument as in the proof of Lemma~\ref{im_quad}.  The set $\{ t+2 : t \in \zed,\ 2 \nmid t \}$ contains the set $\pee$ of all odd prime numbers. Then select $t \geq 3$ such that $p=t+2 \in \pee$, so
$$D=p(p-4),$$
and clearly $\gcd(p,p-4)=1$. To ensure that $D$ is squarefree we need to select $t$ such that $p-4$ is squarefree. The fact that there exist infinitely many prime numbers $p$ such that $p-4$ is squarefree again follows from the theorem on p. 920 of \cite{clary}; for each such prime $p$, let $t=p-2$. For each such choice of $t$, let $K=\que(\sqrt{D})$ and 
$$I = \left< 2b+1, b+\delta \right>  = \left< t+2, \frac{t+2}{2} - \frac{\sqrt{D}}{2} \right> \subseteq \O_K$$
be an ideal. Then
\begin{eqnarray*}
\NN(b+\delta) & = & \left(\frac{(t+2)-\sqrt{D}}{2} \right) \left(\frac{(t+2)+\sqrt{D}}{2} \right) = \frac{1}{4} (t^2 + 4t + 4 - D) \\
& = & t+2 = a,
\end{eqnarray*}
and so condition \eqref{norm_abg} is satisfied. Therefore $t+2, \frac{t+2}{2} - \frac{\sqrt{D}}{2}$ is a canonical basis for $I$. Then 
$$\Lambda_K(I) = \left( \begin{matrix} t+2 & \frac{t+2}{2} - \frac{\sqrt{D}}{2} \\ t+2 & \frac{t+2}{2} + \frac{\sqrt{D}}{2} \end{matrix} \right) \zed^2 = A\zed^2,$$
where we make a change of basis so that
\begin{equation}
\label{A_2}
A = \left( \begin{matrix} \frac{t+2}{2} + \frac{\sqrt{D}}{2} & \frac{t+2}{2} - \frac{\sqrt{D}}{2} \\ \frac{t+2}{2} - \frac{\sqrt{D}}{2} & \frac{t+2}{2} + \frac{\sqrt{D}}{2} \end{matrix} \right).
\end{equation}
Then for any $\bx = A \left( \begin{matrix} m \\ n \end{matrix} \right) \in \Lambda_K(I)$,
\begin{eqnarray}
\label{Q_mn_r}
\|\bx\|^2 & = & Q(m,n) := (m\ n) A^t A \left( \begin{matrix} m \\ n \end{matrix} \right) \nonumber \\
& = & t(t+2) m^2 + 4(t+2) mn + t(t+2) n^2,
\end{eqnarray}
and hence the positive definite integral binary quadratic form $Q(m,n)$ is reduced and symmetric for each $t \geq 5$. As in the proof of Lemma~\ref{im_quad}, this implies that $\Lambda_K(I)$ is WR. Thus we have constructed an infinite family of real quadratic number fields $\que(\sqrt{D})$ with $D>1$, $D \equiv 1 (\md 4)$, each of which contains at least one ideal $I \subseteq \O_K$ so that $\Lambda_K(I)$ is WR. This completes the proof of the lemma.
\endproof

\begin{rem} \label{basis_change} It is interesting to notice that the basis choices for the lattice $\Lambda_K(I)$ resulting in the reduced symmetric norm form in Lemmas~\ref{im_quad} and \ref{real_quad} are different: in the imaginary quadratic case this basis is \eqref{A_1}, which corresponds to the canonical basis for the ideal, while in the real quadratic case this is the basis \eqref{A_2}, which is obtained from the canonical basis by an elementary change of basis operation.
\end{rem}
\smallskip

\proof[Proof of Theorem \ref{q_ideals}]
Theorem  \ref{q_ideals} now follows immediately from Lemmas~\ref{im_quad} and~\ref{real_quad}.
\endproof
\smallskip

\begin{rem} \label{rem_rq} In Tables \ref{table1} and \ref{table2} we present a few examples of ideals $I$ in quadratic number fields $K=\que(\sqrt{-D})$ with $-D \equiv 1 (\md 4)$ and $K=\que(\sqrt{D})$ with $D \equiv 1 (\md 4)$, respectively, so that $\Lambda_K(I)$ is WR, as discussed in Lemmas~\ref{im_quad} and \ref{real_quad}; for each such ideal, we give a presentation in terms of the canonical basis and explicitly write down elements of $I$ that result in minimal vectors in $\Lambda_K(I)$ under the embedding~$\sigma$ (we call them {\it minimal elements}). Notice that these families consists only of some ideals for which the quadratic form $Q$, either corresponding to the unique choice of the basis as in \eqref{I_abg} or obtained from it by one elementary change of basis operation, is reduced and symmetric. There may of course be many other such examples, as well as other more complicated situations when the form is not reduced, but is equivalent to a symmetric reduced form, in which case the lattice in question is again WR. In other words, there are likely many more WR lattices coming from ideals in real and imaginary quadratic fields than the proofs of Lemmas \ref{im_quad} and \ref{real_quad} demonstrate. In order for this to happen, it is necessary for the discriminant of the corresponding positive definite quadratic form $Q$ to have a class represented by a symmetric reduced form; see pp. 19--20 of \cite{buell} for some computational data and p. 27 for related remarks.
\end{rem}
\smallskip

\begin{center} 
\begin{table}[!h!b!p!]
\caption{Examples of ideals in imaginary quadratic fields $K=\que(\sqrt{-D})$ that give rise to WR lattices.} 
\begin{tabular}{|l|l|l|l|} \hline
{\em $-D$} & {Ideal $I \subset \O_K$} & Minimal elements \\ \hline \hline
-15 & $\left< 2,\frac{1-\sqrt{-15}}{2} \right>$ & $\pm 2, \pm \frac{1-\sqrt{-15}}{2}$ \\ \hline
-55 & $\left< 4,\frac{3-\sqrt{-55}}{2} \right>$ & $\pm 4, \pm \frac{3-\sqrt{-55}}{2}$ \\ \hline
-119 & $\left< 6,\frac{5-\sqrt{119}}{2} \right>$ & $\pm 6, \pm \frac{5-\sqrt{119}}{2}$ \\ \hline
-207 & $\left< 8,\frac{7-\sqrt{207}}{2} \right>$ & $\pm 8, \pm \frac{7-\sqrt{207}}{2}$ \\ \hline
\end{tabular}
\label{table1}
\end{table}
\end{center}

\begin{center} 
\begin{table}[!ht]
\caption{Examples of ideals in real quadratic fields $K=\que(\sqrt{D})$ that give rise to WR lattices.} 
\begin{tabular}{|l|l|l|l|} \hline
{\em $D$} & {Ideal $I \subset \O_K$} & Minimal elements \\ \hline \hline
21 & $\left< 7,\frac{7-\sqrt{21}}{2} \right>$ & $\pm \frac{7 \pm \sqrt{21}}{2}$ \\ \hline
165 & $\left< 15,\frac{15-\sqrt{165}}{2} \right>$ & $\pm \frac{15 \pm \sqrt{165}}{2}$ \\ \hline
285 & $\left< 19,\frac{19-\sqrt{285}}{2} \right>$ & $\pm \frac{19 \pm \sqrt{285}}{2}$ \\ \hline
957 & $\left< 33,\frac{33-\sqrt{957}}{2} \right>$ & $\pm \frac{33 \pm \sqrt{957}}{2}$ \\ \hline
\end{tabular}
\label{table2}
\end{table}
\end{center}

Finally notice that all examples in Tables~\ref{table1} and \ref{table2} have 4 minimal elements. In fact, all elements of the infinite family of ideals we constructed in the proof of Lemma~\ref{im_quad} have 4 minimal elements. We remark on this in our next lemma.

\begin{lem} \label{sqrt3} Let $D \neq \pm 3$ be a squarefree integer and $K = \que(\sqrt{D})$ be a quadratic number field. Let $I \subset \O_K$ be an ideal, then $|S(\Lambda_K(I))| \leq 4$.
\end{lem}

\proof
Let $D$ be a squarefree integer, $K= \que(\sqrt{D})$ be a quadratic number field, and $I \subset \O_K$ an ideal. By Lemma~\ref{pairs}, $|S(\Lambda_K(I))| \leq 4$ unless $\Lambda_K(I) \sim \Lambda_h$, in which case $S(\Lambda_K(I))$ contains 6 vectors. Assume this is the case, then there must exist $\bx,\bwy \in S(\Lambda_K(I))$ such that the angle between these two vectors is $\pi/3$. In other words, one of these two vectors, say $\bwy$, is obtained from the other, $\bx = \left( \begin{matrix} x_{11} + x_{12} \sqrt{|D|} \\ x_{21} + x_{22} \sqrt{|D|} \end{matrix} \right)$ with $x_{ij} \in \que$, by rotating it by $\pi/3$, i.e.
$$\bwy = \left( \begin{matrix} \frac{1}{2} & -\frac{\sqrt{3}}{2} \\ \frac{\sqrt{3}}{2} & \frac{1}{2} \end{matrix} \right) \bx \in \que \left( \sqrt{|D|} \right)^2.$$
This readily implies that $\sqrt{3} \in \que(\sqrt{|D|})$, meaning that $D = \pm 3$. This completes the proof.
\endproof

\begin{rem} \label{hex_sim} Both number fields $\que(\sqrt{-3})$ and $\que(\sqrt{3})$ contain ideals giving rise to WR lattices with six minimal vectors, i.e., similar to the hexagonal lattice. The fact that this is true for every ideal of $\que(\sqrt{-3})$ is well-known (in particular, it follows from our Corollary~\ref{cycl_pid}). As for $\que(\sqrt{3})$, it is easy to see that if $I = \left< 1-\sqrt{3} \right> \subset \O_{\que(\sqrt{3})}$, then
$$\Lambda_K(I) = \left( \begin{matrix} 1+\sqrt{3} & 2 \\ 1-\sqrt{3} & 2 \end{matrix} \right) \zed^2$$
is a WR lattice with $S(\Lambda_K(I)) = \left\{ \pm \left( \begin{matrix} 2 \\ 2 \end{matrix} \right), \pm \left( \begin{matrix} 1+\sqrt{3} \\ 1-\sqrt{3} \end{matrix} \right), \pm \left( \begin{matrix} 1-\sqrt{3} \\ 1+\sqrt{3} \end{matrix} \right) \right\}$. This lattice is similar to $\Lambda_h$.
\end{rem}
\bigskip

\section{WR principal ideal lattices}
\label{cyclotom}

In this section we investigate principal ideal lattices with the WR property in any dimension, proving that they come only from cyclotomic fields. We first need a simple minimization lemma, which is essentially the extremal case of the arithmetic mean - geometric mean inequality (from now on abbreviated AM-GM inequality).

\begin{lem} \label{optimize} Let $N \geq 1$ be an integer, and write $\bY=(Y_1,\dots,Y_N)$ for a variable vector. Then
\begin{equation}
\label{amgm}
\left( \prod_{n=1}^N Y_n \right)^{1/N} \leq \frac{1}{N} \sum_{n=1}^N Y_n
\end{equation}
holds for all nonnegative real values of $Y_1,\dots,Y_N$. Moreover, let $A$ be a positive real number, and let
$$f(\bY) = \sum_{n=1}^N Y_n,\ g(\bY) = \prod_{n=1}^N Y_n - A.$$
Then the minimum of $f(\bY)$ under the constraints $g(\bY)=0$ and $Y_n > 0$ for all $1 \leq n \leq N$ is achieved if and only if 
$$Y_1 = \dots = Y_N,$$
in which case
$$\frac{1}{N} \sum_{n=1}^N Y_n = \left( \prod_{n=1}^N Y_n \right)^{1/N} = A^{1/N}.$$
In particular, if $A=1$, this happens when $Y_n=1$ for all $1 \leq n \leq N$, in which case the minimum of the sum $f(\bY)$ is equal to~$N$.
\end{lem}

\proof
Formula \eqref{amgm} is the usual AM-GM inequality. The rest of the statement of this lemma is readily verified, for instance by the method of Lagrange multipliers. This is simply the well-known fact that the arithmetic mean of $N$ positive numbers with a fixed geometric mean is minimized when all of these numbers are equal. The standard geometric interpretation of this fact is that among all $N$-dimensional boxes with a fixed volume, the sum of lengths of edges connected to each vertex is minimized in an $N$-dimensional cube of that volume.
\endproof
\smallskip

In order to proceed with the remainder of this section, we need to set up some basic notation of absolute values on number fields. We always write $K$ for a number field of degree $d$ over $\que$ with $r_1$ real and $2r_2$ complex embeddings and $\O_K$ for its ring of integers, as specified in Section~\ref{intro}.  Let $M(K)$ be the set of places of $K$. For each place $v \in M(K)$ we write $K_v$ for the completion of $K$ at $v$ and let $d_v = [K_v:\que_v]$ be the local degree of $K$ at $v$. Then for each place $u \in M(\que)$ we have
\begin{equation}
\sum_{v \in M(K), v|u} d_v = d.
\end{equation}
For each place $v \in M(K)$ we define the absolute value $|\ |_v$ to be the unique absolute value on $K_v$ that extends either the usual absolute value $|\ |$ on $\real$ or $\cee$ if $v | \infty$, or the usual $p$-adic absolute value on $\que_p$ if $v|p$, where $p$ is a prime. Therefore, the archimedean places are split into the real $v_1,\dots,v_{r_1}$ and the complex $u_1,\dots,u_{r_2}$, given by
$$|x|_{v_n} = |\sigma_n(x)|,\ \forall\ 1 \leq n \leq r_1,$$
and
$$|x|_{u_m} = |\tau_m(x)| = |\bar{\tau}_m(x)| = \sqrt{ \tau_{m1}(x)^2 + \tau_{m2}(x)^2 },\ \forall\ 1 \leq m \leq r_2,$$
for each $x \in K$. For every finite place $v \in M(K)$, $v \nmid \infty$, we define the {\it local ring of $v$-adic integers} $\OO_v = \{ x \in K : |x|_v \leq 1 \}$, whose unique maximal ideal is $\MM_v =  \{ x \in K : |x|_v < 1 \}$. Then $O_K = \bigcap_{v \nmid \infty} \OO_v$. For each $0 \neq x \in K$ the {\it product formula} reads
\begin{equation}
\label{product_formula}
\prod_{v \in M(K)} |x|^{d_v}_v = 1.
\end{equation}
Then for each $0 \neq x \in \O_K$ we must have
\begin{equation}
\label{norm_x}
|\NN(x)| = \prod_{v \mid \infty} |x|_v^{d_v} \geq 1, \text{ since } \prod_{v \nmid \infty} |x|_v^{d_v} \leq 1,
\end{equation}
where $\NN(x)$ is the norm of $x$.
\smallskip

With this notation at hand, we can proceed to our next lemma, which provides a basic description of the set of minimal vectors for a principal ideal lattice.

\begin{lem} \label{lattice_min} Let $K$ be a number field of degree $d$ over $\que$ with $r_1$ real and $2r_2$ complex embeddings, and let $I \subseteq \O_K$ be an ideal. Then
\begin{equation}
\label{min_norm_ideal}
|\Lambda_K(I)| \geq (r_1+r_2) \NN(I)^{\frac{1}{r_1+r_2}},
\end{equation}
and $|\Lambda_K| = r_1+r_2$. Moreover, if $\sigma(x) \in S(\Lambda_K)$ for some $x \in \O_K$, then $x$ is a root of unity. Hence we have
$$S(\Lambda_K) = \{ \bx \in \Lambda_K : \|\bx\|^2 = r_1+r_2 \} = \{ \sigma(x) : x \in \O_K \text{ is a root of unity} \}.$$
\end{lem}

\proof Fix $0 \neq x \in \O_K$, and order the real embeddings in such a way that
\begin{equation}
\label{sigmas}
|\sigma_1(x)|,\dots,|\sigma_k(x)| \geq 1,\ |\sigma_{k+1}(x)|,\dots,|\sigma_{r_1}(x)| < 1,
\end{equation}
for some $0 \leq k \leq r_1$. Then
\begin{eqnarray}
\label{abs_val_ineq}
&\ & 1\ \leq\  |\NN(x)|^{\frac{1}{r_1+r_2}} = \left\{ \prod_{v \mid \infty} |x|_v^{d_v} \right\}^{\frac{1}{r_1+r_2}} \nonumber \\
& \ &\ \ \ = \ \left( \prod_{n=1}^k |\sigma_n(x)| \times \prod_{n=k+1}^{r_1} |\sigma_n(x)| \times \prod_{m=1}^{r_2} (\tau_{m1}(x)^2 + \tau_{m2}(x)^2) \right)^{\frac{1}{r_1+r_2}} \\
&\ &\ \ \ \leq \left( \prod_{n=1}^k |\sigma_n(x)|^2 \times \prod_{n=k+1}^{r_1} \left( \frac{1+|\sigma_n(x)|^2}{2} \right) \times \prod_{m=1}^{r_2} (\tau_{m1}(x)^2 + \tau_{m2}(x)^2) \right)^{\frac{1}{r_1+r_2}} \nonumber
\end{eqnarray}
where the last inequality follows by the AM-GM inequality, since for every real number $a \geq 0$,
$$a = \sqrt{1 \times a^2} \leq \frac{1+a^2}{2}.$$
Applying the AM-GM inequality to \eqref{abs_val_ineq} once again, we see that $|\NN(x)|^{\frac{1}{r_1+r_2}}$ is
\begin{eqnarray}
\label{ar_geom_ineq}
&\ & \leq\ \frac{1}{r_1+r_2} \left( \sum_{n=1}^{k} \sigma_n(x)^2 + \frac{1}{2} \sum_{n=k+1}^{r_1} (1+\sigma_n(x)^2) + \sum_{m=1}^{r_2} (\tau_{m1}(x)^2 + \tau_{m2}(x)^2) \right) \nonumber \\
&\ & =\ \frac{1}{r_1+r_2} \left( \|\sigma(x)\|^2 + \frac{1}{2} \sum_{n=k+1}^{r_1} (1 - \sigma_n(x)^2) \right) \leq  \frac{\|\sigma(x) \|^2}{r_1+r_2},
\end{eqnarray}
which implies that 
\begin{equation}
\label{sr1r2}
\|\sigma(x)\|^2 \geq (r_1+r_2) |\NN(x)|^{\frac{1}{r_1+r_2}} \geq r_1+r_2
\end{equation}
for each $x \in \O_K$. On the other hand, one readily checks that $\|\sigma(1)\|^2 = r_1+r_2$, and hence $|\Lambda_K|=r_1+r_2$. Hence if $\sigma(x) \in S(\Lambda_K)$, we must have $|\NN(x)|=1$, meaning that $x$ is a unit. Now suppose that $I \subseteq \O_K$ is an ideal and $x \in I$. By Lemma~5.1 of \cite{cohn}, we see that
\begin{equation}
\label{norm_ineq}
|\NN(x)| \geq \NN(I).
\end{equation}
Then \eqref{sr1r2} implies that $\|\sigma(x)\|^2 \geq (r_1+r_2) \NN(I)^{\frac{1}{r_1+r_2}}$, which proves \eqref{min_norm_ideal}. 

For the rest of this proof, assume that $x \in \O_K$ is such that $\sigma(x) \in S(\Lambda_K)$, so $\|\sigma(x)\|^2 = r_1+r_2$, which is the smallest possible. Then, combining \eqref{abs_val_ineq} with \eqref{ar_geom_ineq}, we obtain
\begin{eqnarray}
\label{main_ineq}
1 & \leq & \left( \prod_{n=1}^k \sigma_n(x)^2 \times \prod_{n=k+1}^{r_1} \left( \frac{1+\sigma_n(x)^2}{2} \right) \times \prod_{m=1}^{r_2} (\tau_{m1}(x)^2 + \tau_{m2}(x)^2) \right)^{\frac{1}{r_1+r_2}} \nonumber \\
& \leq & \frac{1}{r_1+r_2} \left( \sum_{n=1}^{k} \sigma_n(x)^2 + \frac{1}{2} \sum_{n=k+1}^{r_1} (1+\sigma_n(x)^2) + \sum_{m=1}^{r_2} (\tau_{m1}(x)^2 + \tau_{m2}(x)^2) \right) \nonumber \\
&\ \ \leq & \frac{\|\sigma(x)\|^2}{r_1+r_2} = 1,
\end{eqnarray}
and hence there must be equality throughout \eqref{main_ineq}. By Lemma~\ref{optimize}, this happens if and only if
\begin{eqnarray}
\label{equal}
& \ & \sigma_1(x)^2 = \dots = \sigma_k(x)^2 \nonumber \\
& \ & = \left( \frac{1+\sigma_{k+1}(x)^2}{2} \right) = \dots = \left( \frac{1+\sigma_{r_1}(x)^2}{2} \right) \nonumber \\
& \ & = \tau_{11}(x)^2 + \tau_{12}(x)^2 = \dots = \tau_{r_21}(x)^2 + \tau_{r_22}(x)^2 = 1.
\end{eqnarray}
Combining this observation with \eqref{sigmas}, we conclude that $k=r_1$, and therefore $x$ must be a root of unity, by Kronecker's Theorem. Conversely, if $x$ is a root of unity, then
$$\sigma_1(x)^2 = \dots = \sigma_{r_1}(x)^2 = \tau_{11}(x)^2 + \tau_{12}(x)^2 = \dots = \tau_{r_21}(x)^2 + \tau_{r_22}(x)^2 = 1,$$
and so $\|\sigma(x)\|^2 = r_1+r_2$, meaning that $\sigma(x) \in S(\Lambda_K)$. This completes the proof of the lemma.
\endproof

Next we show that principal ideal lattices coming from cyclotomic fields are always WR.



\begin{lem} \label{cyclotomic} Let $k$ be a positive integer, let $\zeta_k$ be primitive $k$-th root of unity, and let $K=\que(\zeta_k)$ be $k$-th cyclotomic field. Then the lattice $\Lambda_K$ is WR in $\real^d$, where $d=\varphi(k)$.
\end{lem}

\proof
If $k=1,2$, then $K=\que$, and so $\Lambda_K=\zed$, which is WR. If $k=3,4$, the result follows from Lemma~\ref{O_K_WR}. If $k > 4$, it is clear that for $K=\que(\zeta_k)$, $r_1=0$, and for each $1 \leq m \leq r_2=d/2$ and $n \in \zed_{\geq 0}$,
$$|\tau_m(\zeta_k^n)|^2 = |\bar{\tau}_m(\zeta_k^n)|^2 = \tau_{m1}(\zeta_k^n)^2+\tau_{m2}(\zeta_k^n)^2 = 1,$$
meaning that $\|\sigma(\zeta_k^n)\|^2 = r_2$, and so $\sigma(\zeta_k^n) \in S(\Lambda_K)$, by Lemma~\ref{lattice_min}. Hence $S(\Lambda_K)$ contains $\varphi(k) = d$ linearly independent vectors, since the collection of elements $\{ \zeta_k^n : 0 \leq n \leq \varphi(k)-1 \}$ forms an integral basis for $\O_K$, and so $\Lambda_K$ is WR.
\endproof

Finally we prove that if $K$ is not cyclotomic, then $\Lambda_K$ cannot be WR.

\begin{lem} \label{only_cyclotomic} Let $K$ be a number field of degree $d=[K:\que] \geq 2$ such that $\Lambda_K$ is WR, then $K$ is cyclotomic.
\end{lem}

\proof
First suppose that $r_1>0$, then the only roots of unity in $K$ are $\pm 1$, meaning that
$$S(\Lambda_K) = \{ \sigma(1), \sigma(-1) \},$$
by Lemma~\ref{lattice_min}. Then $\dim_{\real} \spn S(\Lambda_K) = 1 < d$, and so $\Lambda_K$ is not~WR.

From now on assume that $K$ is totally imaginary, then $r_1=0$ and $2r_2=d=[K:\que]$. Let $x \in \O_K$ such that $\sigma(x) \in S(\Lambda_K)$, then Lemma~\ref{lattice_min} implies that $x$ is a root of unity. Let $\zeta_k$ be a root of unity of the highest order in $K$, then all other roots of unity in $K$ are powers of $\zeta_k$, and so are contained in $\zed[\zeta_k] \subseteq \O_K$, which means that they can be expressed as integral linear combinations of $1,\zeta_k,\dots,\zeta_k^{\varphi(k)-1}$, an integral basis for $\zed[\zeta_k]$. Hence at most $\varphi(k)$ roots of unity in $K$ can be linearly independent over $\zed$, and so $S(\Lambda_K)$ can contain at most $\varphi(k)$ linearly independent vectors. In order for $\Lambda_K$ to be WR, $\varphi(k)$ must be equal to $d$ by Lemma~\ref{lattice_min}, meaning that $K$ is the cyclotomic field $\que(\zeta_k)$.
\endproof
\smallskip

\proof[Proof of Theorem~\ref{main}] 
Theorem~\ref{main} follows upon combining Lemmas~\ref{lattice_min}, \ref{cyclotomic}, and \ref{only_cyclotomic}.
\endproof

As a consequence of Theorem~\ref{main}, we can also prove that in fact all ideal lattices coming from cyclotomic fields are WR.
\smallskip

\proof[Proof of Corollary~\ref{cycl_ideal}]
Let $I$ be a fractional ideal of $\O_K$, where $K=\que(\zeta_k)$ for some primitive $k$-th root of unity $\zeta_k$, $k \geq 2$. If $k=2$, the result follows from Corollary \ref{cycl_pid}, so assume that $k \geq 3$. Suppose that $\sigma(x) \in S(\Lambda_K(I))$ for some $x \in I$. Notice that, since $K$ is a cyclotomic field,
$$\|\sigma(x)\|^2 = \sum_{n=1}^{\varphi(k)} \tau_n(x) \bar{\tau}_n(x).$$
Then for each $0 \leq m \leq k$,
\begin{eqnarray}
\label{ru_equal}
\|\sigma(\zeta_k^m x)\|^2 & = & \sum_{n=1}^{\varphi(k)} \tau_n(\zeta_k^m x) \bar{\tau}_n(\zeta_k^m x) \nonumber \\
& = & \sum_{n=1}^{\varphi(k)} \tau_n(\zeta_k^m) \bar{\tau}_n(\zeta_k^m) \tau_n(x) \bar{\tau}_n(x) \nonumber \\
& = & \sum_{n=1}^{\varphi(k)} \tau_n(x) \bar{\tau}_n(x) = \|\sigma(x)\|^2,
\end{eqnarray}
and hence $\sigma(\zeta_k^m x) \in S(\Lambda_K(I))$ for each $0 \leq m \leq k$. Since the collection of elements $\{ \zeta_k^m : 0 \leq m \leq \varphi(k)-1 \}$ forms an integral basis for $\O_K$, the collection $\{ \zeta_k^m x : 0 \leq m \leq \varphi(k)-1 \}$ is linearly independent, which in turn implies that the collection of vectors
$$\{ \sigma(\zeta_k^m x) : 0 \leq m \leq \varphi(k)-1 \} \subset S(\Lambda_K(I))$$
is linearly independent in $\real^{\varphi(k)}$, and so $\Lambda_K(I)$ is WR.
\endproof
\bigskip

{\bf Acknowledgment.} We would like to thank Professors Wai Kiu Chan and Sinnou David for their highly helpful comments on the subject of this paper. We would also like to thank the referee for the very useful remarks.
\bigskip

\bigskip

\bibliographystyle{plain}  
\bibliography{wr_ideal}    

\begin{thebibliography}{10}

\bibitem{esm}
A.~H. Banihashemi and A.~K. Khandani.
\newblock On the complexity of decoding lattices using the {K}orkin-{Z}olotarev
  reduced basis.
\newblock {\em IEEE Trans. Inform. Theory}, 44(1):162--171, 1998.

\bibitem{bayer1}
E.~Bayer-Fluckiger.
\newblock Lattices and number fields.
\newblock Contemp. Math. 241, pages 69--84, 1999.

\bibitem{bayer2}
E.~Bayer-Fluckiger.
\newblock Ideal lattices.
\newblock In {\em A panorama of number theory or the view from Baker's garden
  (Zurich, 1999)}, pages 168--184. Cambridge Univ. Press, Cambridge, 2002.

\bibitem{bayer_eucl}
E.~Bayer-Fluckiger.
\newblock Upper bounds for {E}uclidean minima of algebraic number fields.
\newblock {\em J. Number Theory}, 121(2):305--323, 2006.

\bibitem{bayer_nebe}
E.~Bayer-Fluckiger and G.~Nebe.
\newblock On the {E}uclidean minimum of some real number fields.
\newblock {\em J. ThŽor. Nombres Bordeaux}, 17(2):437Ð454, 2005.

\bibitem{oggier}
E.~Bayer-Fluckiger, F.~Oggier, and E.~Viterbo.
\newblock Algebraic lattice constellations: bounds on performance.
\newblock {\em IEEE Trans. Inform. Theory}, 52(1):319--327, 2006.

\bibitem{buell}
D.~A. Buell.
\newblock {\em Binary Quadratic Forms}.
\newblock Springer-Verlag, 1989.

\bibitem{clary}
S.~Clary and J.~Fabrykowski.
\newblock Arithmetic progressions, prime numbers, and squarefree integers.
\newblock {\em Czechoslovak Math. J.}, 54(129)(4):915--927, 2004.

\bibitem{cohn}
H.~Cohn and N.~Heninger.
\newblock Ideal forms of {C}oppersmith's theorem and {G}uruswami-{S}udan list
  decoding.
\newblock {\em preprint, arXiv:math.NT/1008.1284v1}.

\bibitem{wr1}
L.~Fukshansky.
\newblock On distribution of well-rounded sublattices of {$\zed^2$}.
\newblock {\em J. Number Theory}, 128(8):2359--2393, 2008.

\bibitem{wr2}
L.~Fukshansky.
\newblock On similarity classes of well-rounded sublattices of {$\zed^2$}.
\newblock {\em J. Number Theory}, 129(10):2530--2556, 2009.

\bibitem{wr3}
L.~Fukshansky, D.~Moore, R.~A. Ohana, and W.~Zeldow.
\newblock On well-rounded sublattices of the hexagonal lattice.
\newblock {\em Discrete Math.}, 310(23):3287--3302, 2010.

\bibitem{marcus}
D.~A. Marcus.
\newblock {\em Number Fields}.
\newblock Springer-Verlag, 1977.

\bibitem{martinet}
J.~Martinet.
\newblock {\em Perfect Lattices in Euclidean Spaces}.
\newblock Springer-Verlag, 2003.

\bibitem{mcmullen}
C.~McMullen.
\newblock Minkowski's conjecture, well-rounded lattices and topological
  dimension.
\newblock {\em J. Amer. Math. Soc.}, 18(3):711--734, 2005.

\bibitem{tsfasman}
M.~A. Tsfasman and S.~G. Vladut.
\newblock {\em Algebraic-Geometric Codes}.
\newblock Kluwer Academic Publishers, 1991.

\end{thebibliography}
\end{document}